\documentclass[10pt]{article}
\usepackage[usenames]{color}
\usepackage{colortbl}
\usepackage{tikz-cd}

\usepackage{xcolor}
\usepackage{amscd,amsmath, amssymb, fancyhdr, mathbbol, dutchcal, url, euscript, amssymb, mathrsfs}
\usepackage[backref=page]{hyperref}
\usepackage{mathtools}
\renewcommand*{\backrefalt}[4]{%
	\ifcase #1 (Not cited.)%
	\or        (Cited on page~#2.)%
	\else      (Cited on pages~#2.)%
	\fi}

\hypersetup{
	colorlinks   = true,
	citecolor    = magenta,
	linkcolor    =blue,
	urlcolor     =magenta	
}

\numberwithin{equation}{section}
\newcommand{\version}{version 1.0,\ \ October, 2024}

\def\eqref#1{(\ref{#1})}

\newcommand{\arrow}{{\:\longrightarrow\:}}
\newcommand{\Z}{{\Bbb Z}}
\def\C{{\Bbb C}}

\newcommand{\Q}{{\Bbb Q}}
\renewcommand{\H}{{\Bbb H}}

\def\1{\sqrt{-1}\:}

\newcommand{\cntrct}                
{\hspace{2pt}\raisebox{1pt}{\text{$\lrcorner$}}\hspace{2pt}}

\newcommand{\calo}{{\cal O}}

\renewcommand{\tilde}{\widetilde}
\renewcommand{\bar}{\overline}
\renewcommand{\phi}{\varphi}
\renewcommand{\epsilon}{\varepsilon}
\renewcommand{\geq}{\geqslant}
\renewcommand{\leq}{\leqslant}

\newcommand{\End}{\operatorname{End}}
\newcommand{\Tot}{\operatorname{Tot}}
\newcommand{\Id}{\operatorname{Id}}

\newcommand{\Pic}{\operatorname{Pic}}

\newcommand{\Lie}{\operatorname{Lie}}

\newcommand{\Tw}{\operatorname{Tw}}

\newcommand{\GL}{\operatorname{GL}}

\newcommand{\U}{\operatorname{U}}

\newcounter{Mycounter}[section]
\newcounter{lemma}[section]
\renewcommand{\thelemma}{{Lemma \thesection.\arabic{lemma}}}
\newcommand{\lemma}{%
    \setcounter{lemma}{\value{Mycounter}}
    \refstepcounter{lemma}
    \stepcounter{Mycounter}
    {\noindent \bf \thelemma:\ }}

\newcounter{claim}[section]
\renewcommand{\theclaim}{{Claim \thesection.\arabic{claim}}}
\newcommand{\claim}{%
    \setcounter{claim}{\value{Mycounter}}
    \refstepcounter{claim}
    \stepcounter{Mycounter}
    {\noindent \bf \theclaim:\ }}

\newcounter{sublemma}[section]
\newcounter{corollary}[section]
\renewcommand{\thecorollary}{{Corollary \thesection.\arabic{corollary}}}
\newcommand{\corollary}{%
    \setcounter{corollary}{\value{Mycounter}}
    \refstepcounter{corollary}
    \stepcounter{Mycounter}
    {\noindent \bf \thecorollary:\ }}

\newcounter{theorem}[section]
\renewcommand{\thetheorem}{{Theorem \thesection.\arabic{theorem}}}
\newcommand{\theorem}{%
    \setcounter{theorem}{\value{Mycounter}}
    \refstepcounter{theorem}
    \stepcounter{Mycounter}
    {\noindent \bf \thetheorem:\ }}

\newcounter{conjecture}[section]
\newcounter{proposition}[section]
\newcounter{definition}[section]
\renewcommand{\thedefinition}
      {{Definition~\thesection.\arabic{definition}}}
\newcommand{\definition}{%
    \setcounter{definition}{\value{Mycounter}}
    \refstepcounter{definition}
    \stepcounter{Mycounter}
    {\noindent \bf \thedefinition:\ }}

\newcounter{example}[section]
\renewcommand{\theexample}{{Example \thesection.\arabic{example}}}
\newcommand{\example}{%
    \setcounter{example}{\value{Mycounter}}
    \refstepcounter{example}
    \stepcounter{Mycounter}
    {\noindent \bf \theexample:\ }}

\newcounter{remark}[section]
\renewcommand{\theremark}{{Remark \thesection.\arabic{remark}}}
\newcommand{\remark}{%
    \setcounter{remark}{\value{Mycounter}}
    \refstepcounter{remark}
    \stepcounter{Mycounter}
    {\noindent \bf \theremark:\ }}

\newcounter{problem}[section]
\newcounter{question}[section]
\newcommand{\proof}{\noindent{\bf Proof:\ }}

\makeatletter

\renewcommand{\leftmark}{{\scriptsize Twistor space of a hypercomplex manifold is never Moishezon}}

\@addtoreset{equation}{section} \@addtoreset{footnote}{section}
\makeatother

\def\blacksquare{\hbox{\vrule width 5pt height 5pt depth 0pt}}
\def\endproof{\blacksquare}

\title{The twistor space of a compact hypercomplex manifold is never Moishezon}
\author{Yulia Gorginyan}
\begin{document}
\maketitle
\begin{abstract} 
 Let $(X,I,J,K)$ be a compact hypercomplex manifold, i.e. a smooth manifold $X$ with an action of the quaternion algebra $\langle Id,I,J,K,\rangle=\H$ on the tangent bundle $TX$, inducing integrable almost complex structures. For any $(a, b, c) \in S^2$, the linear combination $L := aI + bJ + cK$ defines another complex structure on $X$. This results in a $\C P^1$-family of complex structures called {\bf the twistor family}.  Its total space is called {\it the twistor space}. We show that the twistor space of a compact hypercomplex manifold is never Moishezon and, moreover, it is never Fujiki class $\mathcal{C}$ (in particular, never K\"ahler and never projective).
\end{abstract}
\tableofcontents
\section{Introduction}
\subsection{4-dimensional world}

The concept of a twistor space was introduced at the end of the 1960s within the realm of theoretical physics, by R. Penrose \cite{Pen}. In differential geometry the framework for these ideas was developed by M. Atiyah, N. Hitchin, and I. Singer in their seminal paper \cite{AHS}. From this perspective, twistor spaces are certain complex $3$-manifolds associated with anti-selfdual (ASD) Riemannian 4-manifolds.

\hfill 

In \cite{Hit} N. Hitchin has shown that there are only two examples of twistor space of compact ASD Riemannian 4-manifolds admitting a K\"ahler structure. These are $\C P^3$, associated with the standard 4-dimensional sphere $S^4$, and the flag space $\mathbb{F}_3(\C)$ associated with the complex projective plane $\C P^2$. 

\hfill

\definition\label{Def_1_Moishezo}
A compact complex manifold $X$ is called {\bf Moishezon} if it is bimeromorphically equivalent to a projective manifold, i.e. there exists a projective manifold $\tilde{X}$ and a holomorphic bimeromorphic map $\mu:\tilde{X}\arrow X$.

\hfill

An example of {\it a non-projective} Moishezon twistor spaces was first produced by Y.S. Poon \cite{Po}. The example was provided by a manifold diffeomorphic to $\C P^2\#\C P^2$ with the corresponding twistor space $Z$ bimeromorphic to $\C P^3$.

\hfill

In the work \cite{Camp1} F. Campana has shown that a twistor space is  Moishezon only when the $4$-manifold is $S^4$ or $\#_n\C P^2$.

\hfill

In fact, the Fubini-Study metric on $\C P^2$ is the unique self-dual structure on $\C P^2$ whose twistor space is Moishezon \cite{Po}. So it  makes sense to study \textquotedblleft good\textquotedblright metrics and ask whether these metrics lead to a Moishezon twistor space. For example, one can study {\it the pluriclosed} metrics. 

\hfill

\definition\label{skt}
A Hermitian metric on a complex manifold is called {\bf SKT} or {\bf pluriclosed} if its fundamental form $\omega$ satisfies $\partial\bar{\partial}\omega=0$.

\hfill

\theorem\label{ration_conned}(Verbitsky, \cite{Ver2})
Let $Z$ be a twistor space of a compact, anti-self-dual Riemannian manifold, admitting an SKT metric. Then $Z$ is K\"ahler, i.e. $\C P^3$ or the flag space in $\C P^2$.

\hfill

The proof of the theorem above is based on the results of \cite{Camp1}. 

\hfill

Later the concept of a twistor space was generalized: the twistor space can be defined for several classes of manifolds with special holonomy.

\hfill

The results in $4$-dimensional geometry often have parallels in a hypercomplex geometry. For instance, the result \eqref{ration_conned} about SKT metrics in the $4$-dimensional case was generalized in a recent work by A. Pipitone Federico \cite{PF} for the hypercomplex setting. A. Pipitone Federico has shown that the twistor space of a compact hypercomplex manifold is never SKT. 

\hfill

In the present paper we focus on twistor spaces of hypercomplex manifolds. 

\subsection{Hypercomplex world}

Recall that {\bf a hypercomplex manifold} is a smooth manifold $M$ equipped with almost complex structures $I,J,K\in\End(TM)$ that are integrable and satisfy the quaternionic relations 
\begin{equation}\label{quaternionic_relations}
    I^2=J^2=K^2=-\Id,\qquad IJ=-JI=K.
\end{equation}
Thus, there is an action of the quaternion algebra on the tangent bundle $TM$.

\hfill

For any $(a, b, c) \in S^2$, the linear combination $L := aI + bJ + cK$ defines another complex structure on $M$.
These linear combinations define {\it a $\C P^1$-family} of complex structures, which is called {\bf the twistor family}.

\hfill

\definition\label{twistor_space_hc}
Let $M$ be a hypercomplex manifold.
{\bf The twistor space} of $M$ is a new complex manifold $\Tw(M)$ diffeomorphic to $ M\times \mathbb{C}\rm{P}^1$ with an almost complex structure defined as follows.
 For any point $(x,L)\in M\times \mathbb{C}\rm{P}^1$ the complex structure on $T_{(x,L)}\Tw(M)$ is given by $L$ on $T_xM$ and the standard complex structure $I_{\mathbb{C}P^1}$ on $T_L\mathbb{C}\rm{P}^1$. This almost complex structure is always integrable \cite{Ob}, \cite{Sal}, \cite{K}. 

 \hfill
 
 The space $\Tw(M)$ is equipped with two projections
 \[\begin{tikzcd}
	& \Tw(M) \\
	M && \C P^1
	\arrow["\sigma",from=1-2, to=2-1]
	\arrow[from=1-2, to=2-3, "\pi", labels=below left],
\end{tikzcd}\]
 where $\pi:\Tw(M)\arrow\mathbb{C}\rm{P}^1$ is holomorphic. The fiber $\pi^{-1}(L)$ at a point $L\in\mathbb{C}\rm{P}^1$ is biholomorphic to the complex manifold $(M,L)$. The fibers of the projection $\sigma$ $\sigma^{-1}(m)=: S_m\cong\C P^1$ are complex submanifolds in $\Tw(M)$ for each $m\in M$; however, the map $\sigma$ is not holomorphic. 

\hfill

We are interested in {\it the algebraic dimension} of a twistor space and list some relevant results and definitions.

\hfill

\definition\label{alg_dim_def} 
Let $X$ be a compact complex manifold. {\bf The algebraic dimension $a(X)$} of $X$
is the transcendence degree over $\C$ of the field of global meromorphic functions $\mathscr{M}(X)$ on $X$.

\hfill

The definition below is equivalent to \ref{Def_1_Moishezo}

\hfill

\definition\label{Def_2_Moishezo} A compact complex manifold $X$ is called {\bf Moishezon} if the algebraic dimension $a(X)$ is equal to the complex dimension $\dim_{\C}(X)$.

\hfill

\definition\label{hk_manifold} A hyperk\"ahler manifold $M$ is a smooth manifold equipped with complex structures $I,J,K\in\End(TM)$ that satisfy the quaternionic relations \eqref{quaternionic_relations} and a Hermitian metric $g$ such that the forms $\Omega_I(X,Y):=g(JX,Y)+\sqrt{-1}g(KX,Y)$ and $\Omega_J(X,Y):=g(KX,Y)+\sqrt{-1}g(IX,Y)$ are closed, $X,Y\in TM$.

\hfill

Intuitively, the Moishezon manifolds are those compact complex manifolds that admit a lot of curves and divisors. These subvarieties can be used to study the geometry of the ambient manifold.
Twistor spaces of compact hyperk\"ahler manifolds are very far from being Moishezon:

\hfill

\theorem\label{alg_dim_twist_hk}(M. Verbitsky, \cite{Ver1})
The twistor space $\Tw(M)$ of a compact hyperk\"ahler manifold $M$ has algebraic dimension $a(\Tw(M))=1$.

\hfill

In the present paper we show that the twistor space of a compact hypercomplex manifold is never Moishezon\footnote{It might be a little counterintuitive since for any given finite set of points in the twistor space there is a rational curve passing through \cite{Ver2}.}. Moreover, it is never of Fujiki class $\mathcal{C}$ \eqref{Fujiki_class_C}.

\hfill

{\bf Acknowledgments:} {I have been profoundly inspired by Misha's discussion on the holography principle in relation to the twistor space of hyperk\"ahler manifolds. He encouraged me to tackle the problem addressed in this paper and provided continuous support throughout its preparation. I would also like to express my gratitude to my friends and colleagues: Nikita Klemyatin, Andrey Soldatenkov, Dmitrii Korshunov, and Bruno Suassuna—for their insightful discussions and invaluable feedback.} 

\section{Preliminaries}
\subsection{Hypercomplex manifolds}

\definition
Let $M$ be a smooth manifold. {\bf An almost complex structure} on $M$ is an endomorphism $I\in\End(TM)$ of the tangent bundle of $M$ satisfying $I^2=-\Id$. The Nijenhuis tensor $N_I$ associated to the almost complex structure $I$ is given by the formula $$N_I(X,Y)=[X,Y]+I[IX,Y]+I[X,IY]-[IX,IY].$$ An almost complex structure is called {\bf integrable} if its Nijenhuis tensor vanishes.


\hfill

\begin{theorem}(Newlander--Nirenberg)
If $I$ is an integrable almost complex structure on $M$, then $M$ admits the structure of a complex manifold compatible with $I$.
\end{theorem}

\hfill

\example\label{Hopf_manifold_1} {\bf A Hopf manifold} $X$ is a compact complex manifold obtained as the quotient of $\C^n\backslash\{0\}$ by a cyclic group $\langle\gamma\rangle$ generated by a holomorphic contraction $\gamma$\footnote{{\bf A holomorphic contraction} with center in $0$ is a map $\gamma:\C^n\arrow\C^n$, $\gamma(0)=0$, such that for some $N\gg 0$ the $N$th iteration $\gamma^N$ maps 
any given compact subset $K\subset\C^n$ of $\C^n$ into an arbitrarily small neighborhood of $0$.}:
\begin{equation*}
    X\cong\C^n\backslash\{0\}/\langle\gamma\rangle
\end{equation*}
 When the dimension $n$ is even and $\gamma\in\GL_n(\H)$, it is, in fact, a compact hypercomplex manifold.  Two-dimensional Hopf manifolds $S$ are called {\bf Hopf surfaces}.

\hfill

The algebraic dimension of the twistor space $Z$ of a Hopf surface $S$ is $a(Z)=2$ \cite{Pon}. We show that the algebraic dimension of the twistor space $\Tw(X)$ of a Hopf manifold $a(\Tw(X))\geq 2$ \eqref{examle_twistors_Hopf_surf}.


\subsection{Hodge structures}

Let $X$ be a compact hypercomplex manifold, $\Tw(X)$ the associated twistor space. We want to study how the cohomology of the fibers of the holomorphic twistor projection $\pi:\Tw(X)\arrow\C P^1$ behaves. To do that, we recall the standard approach to studying of the cohomology in  families.

\hfill

\definition
Fix a natural number $k$. 
{\bf A Hodge structure of weight $k$} is a finitely generated free $\Z$-module $V$ of finite rank together with a direct sum decomposition on $V_{\C}=V_{\Z}\otimes_{\Z}\C$:
\begin{equation}\label{def_hodge_structure}
V_{\C}=\bigoplus_{p+q=k}V^{p,q},\quad\text{with}\quad V^{p,q}=\overline{V^{q,p}}\quad
\end{equation}
The decomposition \ref{def_hodge_structure} is called {\bf the Hodge decomposition}.

\hfill


\hfill

 A Hodge structure is equipped with {\bf a $\U(1)$-action}, with $z\in\U(1)$ acting as $z^{p-q}$ on $V^{p,q}$-subspace.

 \hfill

\definition\label{polarization}  {\bf A polarization} of a Hodge structure $V$ of weight $k$ is a $\U(1)$-invariant bilinear form $Q:V\otimes V\arrow\Z$, such that for its $\C$-bilinear extension to $V_{\C}$
\begin{enumerate}
    \item $Q(u,v)=(-1)^kQ(v,u)$;
    \item $Q(u,v)=0$ for $u\in V^{p,q}, v\in V^{a,b}$, where $p\ne b$ and $q\ne a$;
    \item the form $u\mapsto(\sqrt{-1})^{p-q}Q(u,\bar{u})$ is positive definite on the $V^{p,q}$-component.
\end{enumerate}

\subsection{Variation of Hodge structures}

\definition Let $f:X\arrow B$ be a differentible map, $X$ and $B$ are $C^{\infty}$-manifolds. The map $f$ is called {\bf a submersion} if at each $x\in X$ the differential $Df_x:T_xX\arrow T_{f(x)}B$ is surjective.

\hfill

\definition {\bf The vertical tangent space} $T_{f}X\subset TX$ of a submersion $f:X\arrow B$ is the kernel $\ker(Df)$ of the differential $Df$.

\hfill

\lemma (Ehresmann) Let $f:X\arrow B$ be a holomorphic, submersive, and proper map, $B$ connected. Then all fibers $X_t$ are diffeomorphic.\,\,\,\endproof

\hfill

In particular, it implies that the cohomology of the fibers is locally constant on $B$.

\hfill

Let $X$ be a compact complex manifold, and $B$ a compact complex curve\footnote{Of course, in general, $B$ does not have to be a curve, but any complex manifold. We assume that it $1$-dimensional, since further it is will be $\C P^1$.}. Assume that $B$ parametrizes {\it the deformations} of $X$, i.e. there is a holomorphic submersive map 
\begin{equation}\label{family_def1}
    f:\mathcal{X}\arrow B
\end{equation}
with a fiber $\mathcal{X}_t=X$ for some $t\in B$.
The map $f$ is called {\bf a family} of complex manifolds. We assume that $f$ is proper. Additionally, we assume that all fibers satisfy the $dd^c$-lemma. For example, we may assume that all $X_t$ are K\"ahler manifolds or Moishezon or Fujiki class $\mathcal{C}$. 

\hfill

\remark Note that the cohomology $H^k(X_t,\Z)$ of these fibers forms {\it a local system} on $B$.

\hfill

\definition\label{local_system} Let $B$ be a topological space and $R$ a coefficient ring. {\bf A local system} on $B$ is a locally constant sheaf of $R$-modules. It is uniquely determined by the action of the fundamental group $\pi_1(B)$ on its fibers. In particular, a local system is trivial when $\pi_1(B)=\{e\}$.

\hfill

Let $\underline{\Z}$ be the constant sheaf on $B$. Take its $k$-th higher direct image under the map \eqref{family_def1} and denote it $\mathbb{V}_{\Z}^k:=R^kf_*\underline{\Z}$. Its stalk at point $t\in B$ is isomorphic to the integer cohomology of the fiber:
\begin{equation*}
  (\mathbb{V}^k_{\Z})_t:=H^k(X_t,\Z).  
\end{equation*}
Denote a free $\Z$-module obtained from $\mathbb{V}^k_{\Z}$ by $$\mathbb{V}^k:=\displaystyle\frac{R^kf_*\underline{\Z}}{\text{torsion}}.$$ 
The fiber of $\mathbb{V}^k$ is isomorphic to $\cong\Z^r$ for some $r\in\mathbb{N}$.

\hfill

Let $\mathbb{V}^k_\C:=\mathbb{V}^k\otimes_{\Z}\C$ be the complexification. 
In order to study the variation of Hodge structures, we have to go from the locally constant sheaf $\mathbb{V}^k_\C$ to the associated smooth vector bundle 
\begin{equation}
\mathcal{V}^k:=C^{\infty}_B\otimes_{\C}\mathbb{V}^k_{\C}.
\end{equation}
It comes with a natural flat connection which is called {\it the Gauss-Manin connection}:
\begin{equation}
\nabla:\mathcal{V}^k\arrow\Omega^1_B\otimes_{C^{\infty}_B}\mathcal{V}^k, 
\end{equation}
where $\Omega^1_B$ is a sheaf of $C^{\infty}$ $1$-forms on $B$.

\hfill

\remark
One can always recover the locally constant sheaf $R^kf_*\C\subset\mathcal{V}^k$ as the subsheaf of $\nabla$-flat sections. Recall that the category of local systems over $B$ is equivalent to the category of vector bundles with a flat connection\footnote{It is called \textquotedblleft  the Riemann-Hilbert correspondence\textquotedblright.} \cite[Theorem 2.50]{OV}.

\subsection{The Gauss-Manin connection and monodromy}

Let $f:X\arrow B$ be a proper smooth submersion with fiber $X_t$. Assume that the base $B$ is connected. Consider the local system $R^if_*\C$ with fiber $H^i(X_t,\C)$ (as a bundle of fiberwise closed forms up to fiberwise exact), and tensor it with the sheaf of smooth functions on $B$.

\hfill

\definition \label{hor_vect_field}
Let $f:X\arrow B$ be a smooth submersion.
A vector field $\tilde{\xi}\in TX$ is called {\bf a horizontal lift} of a field $\xi\in TB$ if  for any $x\in X$: $$Df_x(\tilde{\xi}_x)=\xi_{f(x)}.$$

\hfill

Restricted to any path $\gamma:[0,1]\arrow B$, the fibration $f:X\arrow B$ becomes trivial. The Lie derivative along the horizontal vector field $\tilde{\xi}$
maps the fiberwise closed form to a fiberwise closed form. Its action on the cohomology classes is well-defined.
Let $\alpha$ be a fiberwise closed form.

\hfill

{\claim\cite[Chapter 24.2]{OV}
The Gauss-Manin connection $\nabla$ takes the class $[\alpha]\in (R^kf_*\C)_{t\in B}$ to the class $[\Lie_{\tilde{\xi}}\alpha]\in (R^kf_*\C)_{\gamma(t)\in B}$:
\begin{equation*}
    \nabla_{\xi}[\alpha]:=[\Lie_{\tilde{\xi}}\alpha].\,\,\, \endproof
\end{equation*}}

\hfill

  \definition {\bf A real variation of Hodge structure (VHS)} over $B$ 
 is a complex vector bundle $(V,\nabla)$ with a flat connection $\nabla$ equipped with a parallel anti-complex involution and a decomposition ${V}^k=\bigoplus_{p+q=k}V^{p,q}$  which satisfies “Griffiths transversality condition”:
 \begin{equation*}
     \nabla_{\xi^{1,0}}(V^{p,q})\subset V^{p,q}\oplus V^{p-1,q+1},\qquad  \nabla_{\xi^{0,1}}(V^{p,q})\subset V^{p,q}\oplus V^{p+1,q-1},
 \end{equation*}
 where $\xi^{1,0}+\xi^{0,1}\in T^{1,0}_B\oplus T^{0,1}_B$ are the vector fields of types $(1,0)$ and $(0,1)$ respectively. 


\hfill

\definition\label{integral_VHS} A VHS $(V,\nabla)$ is called {\bf integral} if $V$ is equipped with a integer lattice $V_{\Z}$, preserved by $\nabla$.

 \hfill

\definition
{\bf A polarized VHS} is an integral  VHS $(V,\nabla)$, ${V}^k=\bigoplus_{p+q=k}V^{p,q}$ such that $\nabla$
 preserves the polarization and the integer lattice $V_{\Z}$.

\hfill

Let $\mathbb{V}^k$ be a Hodge structure of weight $k$, $t\in B$ a base point.
The homomorphism 
\begin{equation}
    \rho:\pi_1(B,t)\arrow\GL(\mathbb{V}^k_t)
\end{equation}
is called {\bf the monodromy representation}. The image $\Gamma:=\rho(\pi_1(B,t))\subseteq\GL(\mathbb{V}^k_t,\Z)$ is called {\bf the monodromy group}.

\hfill

\theorem (P. Deligne, \cite[Theorem 3.1]{Voi2})\label{Deligne_thm} 
Let $B$ be a smooth quasiprojective variety, $\mathbb{V}$ a polarized VHS on $B$ with the trivial monodromy. Then $\mathbb{V}$ is trivial, i.e. $\mathbb{V}\cong V_0\otimes_{\Z}\underline{\Z}$, where $V_0$ is a fixed Hodge structure.

\subsection{Cohomology of the blow-up}

Let $\mu:\tilde{X}\arrow X$ be a bimeromorphic map. For example, assume that $\mu$ is the blow-up of $X$ along a smooth submanifold $Y$ of complex codimension $r$. Cohomology of $\tilde{X}$  are described in the following way \cite[Theorem 7.31 ]{Voi1}:
\begin{equation*}\label{cohomology_of_the_blowup}
    H^k(\tilde{X},\Z)\cong H^k(X,\Z)\oplus\bigoplus_{i=0}^{r-2}H^{k-2i-2}(Y,\Z)
\end{equation*}

The map of the cohomology $\mu^*$ is an embedding:
\begin{equation*}\label{cohom_of_blow_up_inject}
    \mu^*:H^k(X,\Z)\arrow  H^k(\tilde{X},\Z).
\end{equation*}


\subsection{Fujiki class $\mathcal{C}$ manifolds}


It is known that Moishezon manifolds satisfy $dd^c$-lemma \cite{Moi}.
Another important class of manifolds for which $dd^c$-lemma holds is {\it Fujiki class $\mathcal{C}$} manifolds.

\hfill

\definition\label{Fujiki_class_C}
A complex manifold $X$ is called {\bf Fujiki class $\mathcal{C}$} if it is bimeromorphically equivalent  to a K\"ahler manifold.

\hfill



For the twistor space of a compact hypercomplex manifold, the notion of being Moishezon is equivalent to being of Fujiki class $\mathcal{C}$.

\hfill

\theorem\label{Campana_Fujiki}(Campana, \cite[Corollaire p. 212]{Camp2}) 
Let $X$ be a compact Fujiki class $\mathcal{C}$ manifold, such that for any two points $x,y\in X$ there exists a connected compact subvariety $S$ of dimension $1$, such that $x,y\in S$. Then $X$ is Moishezon.

\hfill

\definition
{\bf An ample rational curve} on a complex manifold is a smooth  curve $C\cong\C P^1$ such that its normal bundle is positive, i.e. $NC\cong\bigoplus\calo(i_k)$, where $i_k>0$.

\hfill

\theorem\label{ample_curve}(\cite[Corollary 4.1]{BC})
The twistor space of a compact hypercomplex manifold contains an ample rational curve.

\hfill

\claim\label{rational_connectedness}(\cite[Claim 2.8]{Ver2})
Let $X$ be a compact complex manifold containing an ample rational curve. Then any $k$ points $z_1,\cdots,z_k$ can be connected by an ample rational curve.

\hfill

\corollary\label{rem_fujiki_hence_moishezon} If the twistor space of a compact hypercomplex manifold is of Fujiki class $\mathcal{C}$, then it is Moishezon.\,\,\, \endproof


\section{Variation of Hodge structures on Moishezon manifolds}


Let $M$ be a Moishezon manifold, $\pi:M\arrow B$ a proper holomorphic submersion. Then the fibers $M_b$ of $\pi$ are also Moishezon manifolds. Consider a local system $\mathbb{V}_{\C}=R^k\pi_*(\C)$.
To show that this construction gives a variation of Hodge structures, it is enough to check the Griffits transversality condition. It is straightforward:

\hfill

\theorem (Griffits transversality condition) Let $\pi:M\arrow B$ be a a proper holomorphic submersion with Moishezon fibers, $\mathbb{V}^k:=R^k\pi_*\C$ the local system, $\nabla:\mathcal{V}^k\arrow\mathcal{V}^k\otimes_{C^{\infty}_B}\Omega^1_B$ the Gauss-Manin connection. Consider the fiberwise Hodge decomposition ${V}^k=\bigoplus_{p+q=k}V^{p,q}$ and let $\xi^{1,0}+\xi^{0,1}\in T^{1,0}_B\oplus T^{0,1}_B$ the vector fields of types $(1,0)$ and $(0,1)$ correspondingly. Then
\begin{equation}
    \nabla_{\xi^{1,0}}(V^{p,q})\subset V^{p,q}\oplus V^{p-1,q+1},\qquad  \nabla_{\xi^{0,1}}(V^{p,q})\subset V^{p,q}\oplus V^{p+1,q-1}.
\end{equation}

\proof
Let $\alpha\in\Omega^{p,q}_{\pi}(M)$ be a fiberwise closed form of type $(p,q)$. The map $\pi:M\arrow B$ is holomorphic, hence it preserves the types. Let $\xi^{1,0}+\xi^{0,1}\in T^{1,0}_B\oplus T^{0,1}_B$ be some vectors on the base, denote the corresponding horizontal lifts by $\xi^{1,0}_{hor},\xi^{0,1}_{hor}\in TM$. Take the class $[\alpha]\in R^{p+q}\pi_*(\C)\otimes C^{\infty}_B$ and apply the Gauss-Manin connection along $\xi^{1,0}$:
\begin{equation}
    \nabla_{\xi^{1,0}}[\alpha]=[\Lie_{\xi_{hor}^{1,0}}\alpha]=[\iota_{\xi_{hor}^{1,0}}d\alpha+d\iota_{\xi_{hor}^{1,0}}\alpha]\in R^{p+q}\pi_*\C
\end{equation}
Note that $d\alpha$ has type $(p+1,q)+(p,q+1)$ and $\iota_{\xi_{hor}^{1,0}}d\alpha$ has type $(p,q)+(p-1,q+1)$. Also, the form $\iota_{\xi_{hor}^{1,0}}\alpha$ has type $(p-1,q)$ and $d\iota_{\xi_{hor}^{1,0}}\alpha$ of type $(p,q)+(p-1,q+1)$. 
\endproof


\section{Hypercomplex Hopf manifolds and their twistor spaces}
 Recall that the twistor space $\Tw(X)$ of a compact hyperk\"ahler manifold $X$ has algebraic dimension $a(\Tw(X))=1$ \eqref{alg_dim_twist_hk}. On the other hand, the twistor space of a compact hypercomplex manifold can have algebraic dimension $>1$, see \ref{examle_twistors_Hopf_surf} below.

\hfill

To provide the example, we start from the following observation.

\hfill

\claim\label{flat_space_twistors}\cite[Claim 2.15]{Ver1}
 Let $\H^n$ be a flat hypercomplex space. Then its twistor space 
\begin{equation*}
    \Tw(\H^n)\cong\Tot\calo(1)^{2n}\cong\C P^{2n+1}\backslash\C P^{2n-1}.\,\,\,\endproof
\end{equation*}

\hfill

\definition\label{alg_reduction} 
Let $X$ be a compact complex manifold. {\bf An algebraic reduction} of $X$ is a projective variety $X^{red}$ together with a meromorphic dominant map $\varphi:X\arrow X^{red}$ such that the associated map $\varphi^*:\mathscr{M}(X^{red})\arrow\mathscr{M}(X)$ of the fields of meromorphic functions is an isomorphism.

\hfill

\remark\label{al_red_triangle}
Let $\varphi:X\arrow X^{red}$ be an algebraic reduction of a compact complex manifold $X$. Then there exists a proper bimeromorphic modification $\mu:\tilde{X}\arrow X$ such that  the diagram commutes: \footnote{the variety $\tilde{X}$ is obtained as the resolution of the closure of the graph of $\varphi$}. 

\[\begin{tikzcd}
	& \tilde{X} \\
	X && X^{red}
	\arrow[from=1-2, to=2-1, "\mu"]
	\arrow[from=1-2, to=2-3,"\psi",  labels=below left]
	\arrow[dashed, from=2-1, to=2-3,"\varphi"],
\end{tikzcd}\]
where the map $\psi:\tilde{X}\arrow X^{red}$ is a proper holomorphic. 

\hfill

\example\label{examle_twistors_Hopf_surf}
Let $X=\H^n\backslash\{0\}/\langle\gamma\rangle$ be a hypercomplex Hopf manifold with the holomorphic contraction acting as 
\[
  \begin{pmatrix}
\lambda \\
 & \ddots\\ 
& & \lambda 
  \end{pmatrix}
\]
for some $\lambda\ne 0$. Such a Hopf manifold admits the natural holomorphic elliptic fibration
\begin{equation}
    X\arrow\C P^{2n-1},\quad (z_1,\cdots,z_{2n})\mapsto [z_1:\cdots z_{2n}]
\end{equation}
with fibers $\C^*/\langle\gamma\rangle$ elliptic curves. 

\hfill

Let $\Tw(\H^n\backslash\{0\})=\Tot\calo(1)^{2n}\backslash\{\text{zero section}\}$ be the twistor space of the punctured flat space $\H^n$, the twistor space of a Hopf manifold $\Tw(\H^n\backslash\{0\}/\langle\gamma\rangle))$. Then we have the following diagram:
\[
\begin{tikzcd}
    & \Tw(\H^n\backslash\{0\}) \\
    \Tw(\H^n\backslash\{0\}/\langle\gamma\rangle) &&\Tw(\H^n\backslash\{0\})/\C^*  \\
    & \C P^1 \times \C P^{2n-1}
    \arrow[from=1-2, to=2-1, "\gamma"]
    \arrow[from=1-2, to=2-3, "\C^*-\text{action}"]
    \arrow[from=2-1, to=3-2, "\C^*/\langle\gamma\rangle-\text{action}"]
    \arrow[from=2-3, to=3-2, "\cong"]
\end{tikzcd}
\]
The bottom left arrow of the diagram above is the holomorphic contraction of the family of elliptic curves; it is a $\C^*/\mu$-principal fiber bundle. 

\hfill

\remark\label{qwe}
We obtain that the algebraic reduction of the space $ \Tw(\H^n\backslash\{0\}/\langle\gamma\rangle)=\Tw(X)$ to $\C P^1 \times \C P^{2n-1}$ (see \ref{al_red_triangle}). 
Note that a meromorphic function on Hopf manifold restricts to a constant on the fiber $\C^*/\langle\gamma\rangle$ (\cite[Proposition 6.5]{Ver4}).
Therefore, the algebraic dimension of the twistor space of an elliptic Hopf manifold is  $a(\Tw(X))=a(\C P^1\times\C P^{2n-1})=2n$.

\hfill

{Note that from \ref{main_theorem} below it follows that  $a(\Tw(X))\leq 2n$. However, the map $\Tw(X)\arrow\C P^1\times\C P^{2n-1}$ is surjective, hence $a(\C P^1\times\C P^{2n-1})\leq a(\Tw(X))$:
$$2n\leq a(\Tw(X))\leq 2n\implies a(\Tw(X))=2n.$$}
This gives an alternative proof of the statement of \ref{qwe}.

\section{The twistor space of a compact hypercomplex manifold}
We will need the following lemma:

\hfill

\lemma\label{Canonical_of_twistors}
Let $X$ be a compact hypercomplex manifold of complex dimension $2n$. The fiberwise canonical bundle $\Omega^{2n}(\Tw(X))$ of the twistor space $\Tw(X)$ is isomorphic to $\calo(-2n)$.

\hfill

\proof
The normal bundle to horizontal section is isomorphic to $\calo(1)^{2n}$, see e.g. \cite[Claim 2.15]{Ver1}.
\endproof

\hfill

\theorem\label{main_theorem}
Let $(X,I,J,K)$ be a compact hypercomplex manifold.
Then $\Tw(X)$ cannot be Moishezon.

\hfill

\proof {\bf Step 1:} We argue ad absurdum. Indeed, assume that $\Tw(X)$ is Moishezon. 
The Hodge-to-de Rham spectral sequence of
Moishezon manifolds degenerates in $E_1$ \cite[Corollary 5.23]{DGMS}. This defines a variation of Hodge structures (VHS)
over $\C P^1$ associated with the Hodge filtration on the cohomology of the fiber of $\pi:\Tw(X)\arrow\C P^1$.

\hfill

{\bf Step 2:} This VHS is in fact polarized. 
Indeed, there exists a bimeromorphic map $\mu:\tilde{\Tw(X)}\arrow \Tw(X)$, such that $\tilde{\Tw(X)}$ projective (\cite{Moi}). 
Hodge structures on projective manifolds admit the natural polarization induced by the Chern class of an ample bundle. By the algebraic version of Sard's lemma the map $$\tilde{\pi}:=\pi\circ\mu:\tilde{\Tw(X)}\arrow\C P^1$$ is a holomorphic submersion over the complement in $\C P^1$ of at most a finite number of points. 
Therefore, the VHS associated with $\pi$ is a polarized subVHS, see \eqref{cohomology_of_the_blowup}.

\hfill

{\bf Step 3:} Deligne's theorem \eqref{Deligne_thm} implies that
any polarized rational VHS over a quasiprojective manifold is trivial if it has trivial monodromy.
Since $\pi_1(\C P^1)$ is a trivial group, the VHS associated with the projection
$\pi:\; \Tw(X) \arrow \C P^1$ is trivial.

\hfill

It remains to prove
that the VHS is non-trivial. 

\hfill

If $H^1(X_I)\neq 0$, we have
a non-trivial VHS on the first cohomology. Indeed, consider two points $\pm I\in \C P^1$ and the corresponding fibers $X_I$ and $X_{-I}$.
The first cohomology of the fibers decompose as follows
\begin{equation}
    H^1(X_I,\Q)\otimes_{\Q}\C=H^{1,0}_I\oplus H^{0,1}_I, \quad H^1(X_{-I},\Q)\otimes_{\Q}\C=H^{1,0}_{-I}\oplus H^{0,1}_{-I}.
\end{equation}
Let $[\alpha]\in H^{1,0}(X_I)$. Then $
I\alpha=\sqrt{-1}\alpha=-(-I)\alpha$. Hence, $-\sqrt{-1}\alpha=(-I)\alpha$ and $[\alpha]\in H^{0,1}(X_{-I})$.
Therefore, the monodromy of VHS is non-trivial.

\hfill

Assume that $H^1(X_I)=0$. In particular, it follows that $H^{0,1}(X_I)=H^1(X_I,\calo_{X_I})=0$. 
Note that the canonical bundle of $X_I$ is topologically trivial. \cite[Section 1.4]{Ver3}.
From the exponential exact sequence of sheaves we get 
\[
\cdots\arrow H^1({X_I},\Z)\arrow H^1({X_I},\calo_{X_I})=0\arrow H^1({X_I},\calo^*_{X_I})\arrow H^2(X,\Z)\arrow\cdots
\]
This implies that
$\Pic^0({X_I})=0$, hence the canonical bundle  of $X_I$ is holomorphically trivial, i.e. there is a holomorphic section $\Phi\in\Omega^{2n}(X_I)$.

\hfill

{\bf Step 4:} Consider the VHS over $\C P^1$ associated
with the middle cohomology of the fiber $X_I$.
By Step 2, it has to be trivial;
by Step 3, $H^{n,0}(X_I)= \langle \Phi\rangle$, and $\Phi$ is a nowhere degenerate holomorphic section of the canonical bundle of $X_I$. In particular, 
\begin{equation*}
    \int_{X_I}\Phi\wedge\bar{\Phi}=\int_{X_{-I}}\Phi\wedge\bar{\Phi}> 0, 
\end{equation*}
i.e. there exists is a non-trivial VHS such that $[\Phi]\in H^{0,n}(X_{-I})$\footnote{Another argument is the following. The fiberwise canonical bundle $\Omega_\pi^{2n}(\Tw(X))$ of $\Tw(X)$ is isomorphic to $\calo(-2n)$ \eqref{Canonical_of_twistors} on each horizontal twistor section, hence the holomorphic $2n$-form $\Phi_I$ cannot be extended to a global fiberwise holomorphic $2n$-form on $\Tw(X)$.}. This implies that the VHS on the middle cohomology of the fibers of $\pi$ is non-trivial, contradicting Step 2. \endproof

\hfill

\corollary The twistor space of a compact hypercomplex manifold is never of Fujiki class $\mathcal{C}$. In particular, it is never K\"ahler.

\hfill

{\proof Follows from \ref{rem_fujiki_hence_moishezon} and \ref{main_theorem}. 
\endproof}









\noindent {\sc Yulia Gorginyan\\
{\sc Instituto Nacional de Matem\'atica Pura e
              Aplicada (IMPA) \\ Estrada Dona Castorina, 110\\
Jardim Bot\^anico, CEP 22460-320\\
Rio de Janeiro, RJ - Brasil\\
\tt  iuliia.gorginian@impa.br \\
}\\

}

\end{document}